\documentclass[12pt]{amsart}
\usepackage{latexsym}
\usepackage{amsfonts}
\newtheorem{prop}{Proposition}
\newcommand{\ovprt}{\overline{\partial}}
\newcommand{\ovli}{\overline}

\begin{document}
\title{The canonical solution operator to $\ovprt $ restricted to 
Bergman spaces}
\author{Friedrich Haslinger}
\date{}
\address{Institut f\"ur Mathematik, Universit\"at Wien,
Strudlhofgasse 4, A-1090 Wien, Austria;
friedrich.haslinger@univie.ac.at \hskip 0.6cm
\newline http://www.mat.univie.ac.at/\~\,has/ }
\subjclass{Primary 32W05; Secondary 32A36}
\keywords{$\ovprt $-equation, Bergman kernel}

\begin{abstract}
We first show that the canonical solution operator to $\ovprt $
restricted  to  $(0,1)$-forms with holomorphic coefficients can
be  expressed by an integral operator using the Bergman kernel.
This  result is used to prove that in the case of the unit disc
in $\mathbb C $ the canonical solution operator to $\ovprt $
restricted   to  $(0,1)$-forms with holomorphic coefficients is
a Hilbert-Schmidt operator.
In  the  sequel  we  give  a direct proof of the last statement
using orthonormal bases and show that in the case of the polydisc
and  the  unit  ball  in   $\mathbb C^n,\  n>1,$ the
corresponding operator fails to be a Hilbert-Schmidt operator.
We also indicate a connection with the theory of Hankel operators.

\end{abstract}
\maketitle

{\bf 1. Introduction.}

Let  $\Omega  $  be  a  bounded domain in $\mathbb C^n$ and let
$A^2(\Omega  )$  denote  the  Bergman  space of all holomorphic
functions $f:\Omega \longrightarrow \mathbb C$ such that
$$\int_{\Omega } |f(z)|^2\,d\lambda(z) < \infty ,$$
where $\lambda $ denotes the Lebesgue measure in $\mathbb C^n .$

We solve the $\overline \partial $-equation $\overline \partial
u  =  g,$  where  $g=\sum_{j=1}^{n}g_j\,d\overline  z_j$  is  a
(0,1)-form  with  coefficents  $g_j  \in  A^2(\Omega ), \ j=1,
\dots n.$

\vskip 0.5 cm

It is pointed out in \cite{FS1} that in the proof that
compactness of the solution
operator for $\ovprt $ on $(0,1)$-forms implies that the boundary of $\Omega $
does not contain any analytic variety of dimension greater than or equal to 1,
only the fact that there is a compact solution operator to $\ovprt $ on the
$(0,1)$-forms with holomorphic coefficients is used.
In this case compactness
of the solution operator restricted to $(0,1)$-forms with holomorphic
coefficients implies already compactness of the solution operator on
general $(0,1)$-forms.

The question of compactness is of interest for various reasons - see 
\cite{FS2} for
an excellent survey.

A similar situation appears in \cite{SSU} where the Toeplitz $C^*$ -algebra
$\mathcal T (\Omega )$ is considered and the relation between the structure of
$\mathcal T (\Omega )$ and the $\ovprt $-Neumann problem is discussed
(see \cite{SSU} , Corollary 4.6).

\vskip 0.3cm

We first show that the canonical solution operator to $\ovprt $
restricted  to  $(0,1)$-forms with holomorphic coefficients can
be  expressed by an integral operator using the Bergman kernel.
This  result is used to prove that in the case of the unit disc
in $\mathbb C ,$ the canonical solution operator to $\ovprt $
restricted   to  $(0,1)$-forms with holomorphic coefficients is
a Hilbert-Schmidt operator.

In  the  sequel  we  give  a direct proof of the last statement
using orthonormal bases and show that in the case of the polydisc
and  the  unit  ball  in  $\mathbb  C^n  \  ,  \  n\ge  2,$ the
corresponding operator fails to be a Hilbert-Schmidt operator.

The canonical solution operator to $\ovprt $
restricted   to  $(0,1)$-forms with holomorphic coefficients
can also be interpreted as the Hankel operator
$$H_{\overline z}(g)=(I-P)(\overline z g),$$
where  $P  :  L^2(\Omega )
\longrightarrow A^2(\Omega )$ denotes the Bergman projection.
See \cite{A}, \cite{AFP}, \cite{B}, \cite{J}, \cite{R},
\cite{W} and \cite{Z} for details.
\vskip 1 cm
{\em Proof.}

{\bf 2. The integral representation.}

The canonical solution operator
$$S_1 : A_{(0,1)}^2(\Omega ) \longrightarrow L^2(\Omega )$$
has  the  properties    $\overline  \partial  S_1  (g) = g$ and
$S_1(g) \perp A^2(\Omega ).$

\begin{prop}
The canonical solution operator
$$S_1 : A_{(0,1)}^2(\Omega ) \longrightarrow L^2(\Omega )$$
has the form
$$S_1(g) (z) = \int_{\Omega }B(z,w)<g(w), z-w>\,d\lambda (w),$$
where $ B$ denotes the Bergman kernel of $\Omega $ and
$$<g(w), z-w>=\sum_{j=1}^n g_j(w) (\overline z_j -\overline w_j
),$$
for $z=(z_1,\dots ,z_n)$ and $w=(w_1,\dots ,w_n).$
\end{prop}
\vskip 0.3 cm

Integral  operators  of  similar  type have been used to settle
questions on compactness of the solution operator to
$\ovprt $, see \cite{CD} and \cite{L}.

\vskip 0.3 cm

{\em Proof.} Let  $v(z)=\sum_{j=1}^n  \overline z_j g_j(z).$ Then it follows
that
$$\overline    \partial   v   =   \sum_{j=1}^n   \frac{\partial
v}{\partial \overline z_j } d\overline z_j =
\sum_{j=1}^n g_j d\overline z_j = g.$$

Hence  the canonical solution operator $S_1 $ can be written in
the  form  $S_1(g)  =  v  -  P(v),$  where  $P  :  L^2(\Omega )
\longrightarrow A^2(\Omega )$ is the Bergman projection.
If $\tilde{v} $ is another solution to
$\overline \partial u = g,$ then $v-\tilde{v} \in A^2(\Omega )$
hence $v=\tilde{v}+h,$ where $h\in A^2(\Omega ).$ Therefore
$$v - P(v) = \tilde{v}+h - P(\tilde{v})-P(h)=
\tilde{v} - P(\tilde{v}).$$

Since $g_j \in A^2(\Omega ), \ j=1,\dots ,n ,$ we have
$$g_j(z)  =  \int_{\Omega  }B(z,w)g_j(w)\,d\lambda (w).$$
Now we get

\begin{eqnarray*}
S_1(g)(z) & =  &\sum_{j=1}^n \overline z_j g_j(z) -
\int_{\Omega }B(z,w)\left (
\sum_{j=1}^n \overline w_j g_j(w) \right ) \,d\lambda (w) \\
& = &  \int_{\Omega }\left [
\left (\sum_{j=1}^n \overline z_j g_j(w) \right ) B(z,w) -
\left (\sum_{j=1}^n \overline w_j g_j(w) \right ) B(z,w)
\right ]\,d\lambda (w)\\
& = & \int_{\Omega }B(z,w)<g(w), z-w>\,d\lambda (w).
\qquad \qquad \Box
\end{eqnarray*}

\vskip 0.3 cm

{\em  Remark.}  It  is pointed  out  that  a  $(0,1)$-form
$g=\sum_{j=1}^{n}g_j\,d\overline  z_j$ with
holomorphic  coefficients  is not invariant under the pull back
by  a  holomorphic  map  $F = (F_1,\dots, F_n):
\Omega_1 \longrightarrow \Omega .$ It can be shown
that
$$F^*g =  \sum_{j=1}^{n}\left ( \sum_{l=1}^n g_l
\frac{\partial \overline F_l}{\partial \overline z_j} \right )
\,d\overline  z_j$$
and the expressions
$\frac{\partial \overline F_l}{\partial \overline z_j}$ are not holomorphic.

Nevertheless it is true that  $\ovprt u = g$ implies
$\ovprt (u \circ F) = F^*g.$

\vskip 1 cm

Now let $ \omega $ be a holomrphic $(n,n)$-form, i.e.
$$\omega = \tilde{\omega }\
d z_1 \wedge \dots \wedge d z_n \wedge d\ovli z_1
\wedge \dots \wedge d\ovli z_n ,$$
where  $\tilde{\omega  }\in  A^2(\Omega ).$ In this case we can
express  the  canonical solution to $\ovprt u = \omega $ in the
following form

\begin{prop}
Let $u$ be the $(n,n-1)$-form
$$u = \sum_{j=1}^n u_j \
d z_1 \wedge \dots \wedge d z_n \wedge d\ovli z_1
\wedge  \dots  \wedge  [d\ovli z_j ] \wedge \dots \wedge d\ovli
z_n,$$
where
$$u_j(z) = \frac{(-1)^{n+j-1}}{n}
\int_{\Omega  }  (\ovli  z_j - \ovli w_j ) B(z,w) \tilde{\omega
}(w)\,d\lambda (w).$$
Then $u_j \perp A^2(\Omega ) \ , j=1,\dots ,n$ and
$\ovprt u = \omega .$
\end{prop}

{\em Proof.} It follows that
$$u_j(z)   =      \frac{(-1)^{n+j-1}}{n}   \left  (  \ovli  z_j
\tilde{\omega  }(z)  -  P(\ovli  w_j \tilde{\omega })(z) \right
),$$
from this we obtain
$$\frac{\partial u_j}{\partial \ovli z_k } =
\frac{(-1)^{n+j-1}}{n}   \left (
\frac{\partial  \ovli  z_j}{\partial  \ovli z_k } \tilde{\omega
}   +  \ovli z_j \frac{\partial \tilde{\omega }}{\partial \ovli z_k}
\right  )  =   \frac{(-1)^{n+j-1}}{n}\ \delta_{jk}\  \tilde{\omega
},$$
where $\delta_{jk}$ is the Kronecker delta symbol. Hence

\begin{eqnarray*}
\ovprt u & = & \sum_{k=1}^n  \sum_{j=1}^n
\frac{\partial u_j}{\partial \ovli z_k }\ d\ovli z_k \wedge
d z_1 \wedge \dots \wedge d z_n \wedge d\ovli z_1
\wedge  \dots  \wedge  [d\ovli z_j ] \wedge \dots \wedge d\ovli
z_n \\
& = & \sum_{k=1}^n  \sum_{j=1}^n
\left ((-1)^{n+j-1}/n \right )\ \delta_{jk}\  \tilde{\omega }\
d\ovli z_k \wedge\\
& \ & \wedge d z_1 \wedge \dots \wedge d z_n \wedge d\ovli z_1
\wedge  \dots  \wedge  [d\ovli z_j ] \wedge \dots \wedge d\ovli
z_n \\
& = & \tilde{\omega }\
d z_1 \wedge \dots \wedge d z_n \wedge d\ovli z_1
\wedge  \dots  \wedge d\ovli z_n. \qquad \qquad \qquad \Box
\end{eqnarray*}

\vskip 0.5 cm
{\em  Remark.}  The  pull  back by a holomorphic map $F$ has in
this case the form
$$F^*\omega = \left |\det \frac{\partial F_j}{\partial z_k}
\right |^2 \ \tilde{\omega }\
d z_1 \wedge \dots \wedge d z_n \wedge d\ovli z_1
\wedge  \dots  \wedge d\ovli z_n.$$

\vskip 0.3 cm
\begin{prop} Suppose that $\Omega $ is a smoothly bounded pseudoconvex
domain of finite type in $\mathbb C^n.$ Let $T: L^2_{(0,1)}(\Omega ) 
\longrightarrow
L^2(\Omega )$ be the operator defined by
$$T(f)(z)= \int_{\Omega }B(z,w)<f(w), z-w>\,d\lambda (w), \ \ f\in
L^2_{(0,1)}(\Omega ).$$
Then $T$ is a compact operator.
\end{prop}
This follows from Theorem 1 in \cite{CD}.

The last result implies that the restriction
of $T$ to $ A^2_{(0,1)}(\Omega ),$ which is the canonical solution operator to
$\ovprt ,$ is also a compact operator. This fact follows also from 
\cite{C}, where
it is shown that the $\ovprt -$ Neumann operator is compact.
\vskip 1 cm

Next we consider the integral kernel of the canonical solution
operator $S_1 $ for the unit disc $\mathbb D $ in $\mathbb C $
and prove that this kernel is square integrable
over $\mathbb D \times \mathbb D .$

\begin{prop}
$$\int_{\mathbb  D}\int_{\mathbb  D}\frac{|\overline z - \overline w
|^2}{|1  -  z\overline  w  |^4}\,d\lambda  (z)\,d\lambda  (w) <
\infty .$$
\end{prop}

{\em Proof.} It is easily seen that
$|z - w|\le |1 - z\overline w| , $ \ for $z,w \in \mathbb D .$
Hence we get
$$\int_{\mathbb  D}\int_{\mathbb  D}\frac{|\overline z - \overline w
|^2}{|1  -  z\overline  w  |^4}\,d\lambda  (z)\,d\lambda  (w)
\le
\int_{\mathbb  D}\int_{\mathbb  D}\frac{1}{|1  -  z\overline
w |^2}\,d\lambda  (z)\,d\lambda  (w) .$$
Using  polar  coordinates  $z  =  r\,e^{i\theta  }$  and  $w  =
s\,e^{i\phi  }$ we can write the last integral in the following
form
$$
\int_{\mathbb  D}\int_{\mathbb  D}\frac{1}{|1  -  z\overline
w |^2}\,d\lambda  (z)\,d\lambda  (w)
  =
\int_0^1\int_0^1\int_0^{2\pi }\int_0^{2\pi }
\frac{r\,s\,d\theta  \,d\phi  \,dr\,ds}{1-2\,r\,s\,\cos (\theta -
\phi ) + r^2\,s^2}$$
$$ =
\int_0^1\int_0^1\int_0^{2\pi }\int_0^{2\pi }
\frac{1-r^2\,s^2}{1-2\,r\,s\,\cos (\theta -
\phi ) + r^2\,s^2}\ \frac{r\,s }{1-r^2\,s^2}
\,d\theta  \,d\phi  \,dr\,ds .$$

Integration of the Poisson kernel with respect to $\theta $ yields
$$\int_0^{2\pi  }\frac{1-\rho^2}{1-2\rho \,\cos(\theta - \phi )
+ \rho^2}\,d\theta = 2\pi \ , 0<\rho <1.$$
Hence
$$\int_0^1\int_0^1\int_0^{2\pi }\int_0^{2\pi }
\frac{1-r^2\,s^2}{1-2\,r\,s\,\cos (\theta -
\phi ) + r^2\,s^2}\ \frac{r\,s }{1-r^2\,s^2}
\,d\theta  \,d\phi  \,dr\,ds $$
$$= (2\pi )^2 \int_0^1\int_0^1 \frac{r\,s }{1-r^2\,s^2}
\,dr\,ds
=  -\,(2{\pi })^2\int_0^1\frac{\log  (1  -  s^2)}{2s}\,ds < \infty .\qquad
\Box $$
\vskip 1 cm

{\em  Remark.}  The  last  proposition implies that the opertor
$T : L^2(\mathbb D ) \longrightarrow L^2(\mathbb D)$ defined by
$$T(f)(z) = \frac{1}{\pi }\,\int_{\mathbb D}\frac{\overline z -
\overline w}{(1-z \overline w )^2}\,f(w)\,d\lambda (w) ,$$
for  $f\in L^2(\mathbb D ),$ is a Hilbert-Schmidt operator, see
\cite{MV} , 16.12.

If   we   restrict   this   operator  to  the  closed  subspace
$A^2(\mathbb D)$ we obtain

\begin{prop} The canonical solution operator to $\ovprt $
$$S_1 : A^2(\mathbb D) \longrightarrow L^2(\mathbb D)$$
is a Hilbert-Schmidt operator.
\end{prop}

{\em  Proof.} By \cite{MV} , 16.8,  we have to show that there exists a
complete  orthonormal  system  $(\phi_k  )_{k=0}^{\infty}$
of $A^2(\mathbb D)$ such that
$$\sum_{k=0}^{\infty }\|S_1(\phi_k )\|^2 < \infty .$$
For this purpose we take a complete orthonormal system
$(\phi_k    )_{k=0}^{\infty}$ of $A^2(\mathbb D)$ and extend it
to a complete orthonormal system $(\psi_j  )_{j=0}^{\infty}$
of    $L^2(\mathbb    D).$  Again  by  \cite{MV}  ,  16.8,  and
proposition 3, it follows that
$$\sum_{j=0}^{\infty }\|T(\psi_j )\|^2 < \infty ,$$
which implies that
$$\sum_{k=0}^{\infty   }\|S_1(\phi_k  )\|^2  <  \infty  .\qquad
\qquad \Box $$

\vskip 1 cm

{\bf 3. Hilbert-Schmidt operators.}
Now we show directly that  the canonical solution operator to
$\ovprt $
$$S_1 : A_{(0,1)}^2(\mathbb D ) \longrightarrow L^2(\mathbb D )$$
is a Hilbert-Schmidt operator, if $\mathbb D$ is the open unit disc in $\mathbb
C,$ and is not Hilbert-Schmidt if $\mathbb B$ is the open unit ball in $\mathbb
C^n $ for $n>1.$
\vskip 0.5 cm
Let $\mathbb D \subset \mathbb C$ and let $\| . \| $ denote the norm in
$A^2(\mathbb D )$
and consider the orthonormal basis
$$\{ u_n(z) = \left [ (n+1)/\pi  \right ]^{1/2} \ z^n \ :
n\in \mathbb N_0 \}$$ of
$A^2(\mathbb D ).$
\begin{prop}
The  canonical  solution  operator  $S_1$  for  the  unit  disc
$\mathbb D $ in $\mathbb C $ has the following property
$$\sum_{n=0}^{\infty }\| S_1 (u_n \ d\ovli z) \|^2 < \infty ,$$
which implies that $S_1 : A_{(0,1)}^2(\mathbb D ) \longrightarrow 
L^2(\mathbb D )$
is a Hilbert-Schmidt operator (see \cite{MV} ).

\end{prop}

{\em Proof.}
Using calculations in \cite{J} we can show that
$$S_1 (u_n \ d\ovli z) (z) = \left [ (n+1)/\pi  \right ]^{1/2}
\  z^n \ \ovli z \ -   \left [ n^2/((n+1)\pi ) \right ]^{1/2} \
z^{n-1}\ , n\in \mathbb N_0.$$

The Bergman kernel $B$ of $\mathbb D $ has the form
$$B(z,  \zeta  ) = \frac{1}{\pi }\ \frac{1}{(1-z\overline \zeta
)^2} ,$$
hence  by  Proposition 1 we can express $\| S_1(u_nd\overline z
)\|^2$ in the form
$$\int_{\mathbb D} \left | \overline z \ u_n(z) -
\frac{1}{\pi  } \int_{\mathbb D}\frac{\overline \zeta \ u_n(\zeta
)}{(1-z\overline    \zeta    )^2}\,d\lambda    (\zeta   )\right
|^2\,d\lambda (z).$$

Therefore we get
\begin{eqnarray*}
\| S_1(u_nd\overline z )\|^2 & = &
\int_{\mathbb  D}\left | \left ( \frac{n+1}{\pi }\right )^{1/2}
\ z^n \ \overline z
- \frac{n\ z^{n-1}}{[(n+1)\pi ]^{1/2}} \right |^2
\,d\lambda (z)\\
& = &
\int_{\mathbb D}\left ( \frac{(n+1)\,|z|^{2n+2}}{\pi }
- \frac{2n\, |z|^{2n}}{\pi } + \frac{n^2\,|z|^{2n-2}}{(n+1)\pi }
\right )\,d\lambda (z)\\
& = &
2 \pi \ \int_0^1
\left ( \frac{(n+1)\, r^{2n+3}}{\pi }
- \frac{2n\, r^{2n+1}}{\pi } + \frac{n^2\, r^{2n-1}}{(n+1)\pi }
\right )\,dr\\
& = &
\frac{1}{(n+1)(n+2)}
\end{eqnarray*}

Hence
$$\sum_{n=0}^{\infty  }\|  S_1(u_n\,d\overline  z) \|^2 < \infty
 \qquad \qquad \Box $$

\vskip 1 cm

{\em Remark.} It can be shown that the set
$\{  S_1(u_n\,d\overline  z) \ : \ n\in \mathbb N_0 \}$ consists
of pairwise orthogonal elements of $L^2 (\mathbb D ).$

\vskip 1 cm

In  the following part we consider the case of the polydisc, in
sake of simplicity we concentrate on $\mathbb C^2 ,$ let
$$\mathbb D^2 = \{ z=(z_1,z_2) \ : \ |z_1|<1 \ , \ |z_2|<1 \}.$$

Now  $\{ z_1^{n_1}\,z_2^{n_2} \ : \ n_1, n_2 \in \mathbb
N_0 \}$ is an orthogonal basis in $A^2(\mathbb D^2 ).$
It is easily seen that the norms of the functions $
z_1^{n_1}\,z_2^{n_2} $ are $\pi [1/((n_1+1)(n_2+1))]^{1/2}.$
The functions
$$u_{n_1,n_2}(z_1,z_2) = \frac{[(n_1+1)(n_2+1)]^{1/2}}{\pi }
\, z_1^{n_1}\,z_2^{n_2} \ , \ n_1, n_2 \in \mathbb N_0$$
form  an  orthonormal  basis  of  $A^2(\mathbb  D^2 ),$ and the
system
$$\{ u_{n_1,n_2}\,d\overline z_1 \ , \
u_{n_1,n_2}\,  d\overline  z_2  \  : \ n_1, n_2 \in \mathbb N_0
\}$$
constitutes an orthonormal basis for
$A^2_{(0,1)}(\mathbb D^2 ).$

Next  we  compute  the  Bergman  projections  of  the  functions
$$(z_1,z_2)\mapsto    \overline    z_1\,   u_{n_1,n_2}(z_1,z_2)
\quad \mbox{and}\quad
(z_1,z_2)\mapsto \overline z_2\, u_{n_1,n_2}(z_1,z_2) $$
and obtain
$$P(\overline \zeta_1\, u_{n_1,n_2}(\zeta_1, \zeta_2))(z_1,
z_2) =
\frac{[(n_1+1)(n_2+1)]^{1/2}}{\pi}\frac{n_1}{n_1+1}\ z_1^{n_1-1}
\,z_2^{n_2}\ ,$$
where we used similar computations as in Proposition 6.

The Bergman projection of the second function is
$$P(\overline \zeta_2 \, u_{n_1,n_2}(\zeta_1, \zeta_2))(z_1,
z_2) =
\frac{[(n_1+1)(n_2+1)]^{1/2}}{\pi}\frac{n_2}{n_2+1}\ z_1^{n_1}
\,z_2^{n_2-1}.$$
Now    we    can  compute  the  norms  of  the images under the
canonical solution operator of the elements of our orthonormal
basis of $A^2_{(0,1)}(\mathbb D^2 ) \,:$
$$\frac{(n_1+1)(n_2+1)}{\pi^2}\, \int_{\mathbb D^2}
\left | \overline z_1\, z_1^{n_1}\, z_2^{n_2} -
\frac{n_1}{n_1+1}\, z_1^{n_1-1}\, z_2^{n_2}\right |^2\,
d\lambda (z)=
\frac{1}{(n_1+2)(n_1+1)},$$
where  we  used  the corresponding computation of
Proposition 6 for the integral
with respect to $z_1.$

In a similar way we obtain
$$\frac{(n_1+1)(n_2+1)}{\pi^2} \int_{\mathbb D^2}
\left | \overline z_2\, z_1^{n_1}\, z_2^{n_2} -
\frac{n_2}{n_2+1}\, z_1^{n_1}\, z_2^{n_2-1} \right |^2\,
d\lambda (z)
  = \frac{1}{(n_2+2)(n_2+1)}.$$
Since
$$\sum_{n_1, n_2 =1}^{\infty}\left (
\frac{1}{(n_1+2)(n_1+1)}+ \frac{1}{(n_2+2)(n_2+1)}\right )
=\infty ,$$
the canonical solution operator
$$S_1  :  A^2_{(0,1)}(\mathbb  D^2 )\longrightarrow L^2(\mathbb
D^2 )$$
is not Hilbert-Schmidt.

{\em Remark.} With results from \cite{K2} it can be shown that
the canonical solution operator
$$S_1  :  A^2_{(0,1)}(\mathbb  D^2 )\longrightarrow L^2(\mathbb
D^2 )$$
is even not compact.
\vskip 1 cm
We  now  consider  the  case  of the unit ball $\mathbb B^2$ in
$\mathbb C^2.$ Here we can use calculations from the proof of Proposition 1
in \cite{W} .

The norms of the functions $
z_1^{n_1}\,z_2^{n_2} $ are now $\pi [n_1!\, n_2!/(n_1+n_2+2)!]^{1/2}$
(see \cite{K1} ).
The functions
$$U_{n_1,n_2}(z_1,z_2) =
\frac{[(n_1+n_2+2)!]^{1/2}}{\pi (n_1!\, n_2!)^{1/2}}
  z_1^{n_1}\,z_2^{n_2} \ , \ n_1, n_2 \in \mathbb N_0$$
form  an  orthonormal  basis  of  $A^2(\mathbb  B^2 ),$ and the
system
$$\{ U_{n_1,n_2}\,d\overline z_1 \ , \
U_{n_1,n_2}\,  d\overline z_2  \  : \ n_1, n_2 \in \mathbb N_0
\}$$
constitutes an orthonormal basis for
$A^2_{(0,1)}(\mathbb B^2 ).$

We compute the Bergman projections of the functions
$$(z_1,z_2)\mapsto    \overline z_1\,  U_{n_1,n_2}(z_1,z_2)
\quad \mbox{and}\quad
(z_1,z_2)\mapsto \overline z_2\, U_{n_1,n_2}(z_1,z_2) $$
and obtain
$$P(\overline \zeta_1\, U_{n_1,n_2}(\zeta_1, \zeta_2))(z_1,
z_2)=
\frac{[(n_1+n_2+2)!]^{1/2}}{\pi\,(n_1!\, n_2!)^{1/2}}
\frac{n_1}{n_1+n_2+2}\,z_1^{n_1-1}z_2^{n_2},$$
$$P(\overline \zeta_2\, U_{n_1,n_2}(\zeta_1, \zeta_2))(z_1,
z_2)=
\frac{[(n_1+n_2+2)!]^{1/2}}{\pi\,(n_1!\, n_2!)^{1/2}}
\frac{n_2}{n_1+n_2+2}\,z_1^{n_1}z_2^{n_2-1}.$$
Finally we  compute  the  norms  of  the images of the basis
elements under the canonical solution opertor $S_1,$ and obtain
$$
\frac{(n_1+n_2+2)!}{\pi^2\,n_1!\, n_2!}
\int_{\mathbb B^2}\left |\overline z_1\,z_1^{n_1}\,z_2^{n_2}
-\frac{n_1}{n_1+n_2+2}\, z_1^{n_1-1}\,z_2^{n_2}\right |^2
\,d\lambda (z)$$
$$ = \frac{n_2+2}{(n_1+n_2+2)(n_1+n_2+3)}$$
and
$$
\frac{(n_1+n_2+2)!}{\pi^2\,n_1!\, n_2!}
\int_{\mathbb B^2}\left |\overline z_2\,z_1^{n_1}\,z_2^{n_2}
-\frac{n_2}{n_1+n_2+2}\, z_1^{n_1}\,z_2^{n_2-1}\right |^2
\,d\lambda (z)$$
$$ =
\frac{n_1+2}{(n_1+n_2+2)(n_1+n_2+3)}.$$

Since
$$\sum_{n_1, n_2 =1}^{\infty}\left (
\frac{n_2+2}{(n_1+n_2+2)(n_1+n_2+3)}+
\frac{n_1+2}{(n_1+n_2+2)(n_1+n_2+3)}\right )
=\infty ,$$
the canonical solution operator
$$S_1  :  A^2_{(0,1)}(\mathbb  B^2 )\longrightarrow L^2(\mathbb
B^2 )$$
is also not Hilbert-Schmidt.
\vskip 0.3 cm
{\em Remark.}In \cite{Z} it is shown that there are no nonzero Hilbert-Schmidt
Hankel operators on the Bergman space of the unit ball in $\mathbb C^n$ with
antiholomorphic symbol when $n\ge 2.$

\vskip 1 cm

\begin{prop} The integral kernel
$$\frac{|z_1-w_1|^2+|z_2-w_2|^2}
{|1-z_1\overline w_1 |^4\,|1-z_2\overline w_2|^4}$$
does not belong to $L^2(\mathbb D^2 \times \mathbb D^2)$ and the 
integral kernel
$$\frac{|z_1-w_1|^2+|z_2-w_2|^2}{|1-z_1\overline w_1 -z_2\overline w_2|^6}
$$
does not belong to $L^2(\mathbb B^2 \times \mathbb B^2).$
\end{prop}

\vskip 1 cm

{\em Proof.} Suppose the first kernel
belongs to $L^2(\mathbb D^2 \times \mathbb D^2)$,
then the  corresponding  integral operator form $L^2_{(0,1)}(\mathbb
D^2  )$  to  $L^2(\mathbb  D^2)$ is a Hilbert-Schmidt operator,
which  would imply that the restriction to $A^2_{(0,1)}(\mathbb
D^2  )$ is also Hilbert-Schmidt, but this restriction coincides
with  the  canonical  solution  operator  $S_1,$  from which we
already  know that it is not Hilbert-Schmidt. The proof for the
second integral is analogous to the first.

\vskip 1 cm

{\em   Acknowledgments.}   The   author   acknowledges   useful
discussions  with  David  Barrett,  Ingo  Lieb and Emil Straube
during a seminar on complex analysis at the Erwin Schr\"o\-dinger
Institute, Vienna.

\end{document}